\documentclass[11pt,reqno]{amsart}

\title[]{On the Space-Time Analyticity of the Inhomogeneous Heat Equation on the Half Space with Neumann Boundary Conditions}

\author{Elie Abdo}
\address{Department of Mathematics, Temple University, Philadelphia, PA 19122, USA}
\email{elie.abdo@temple.edu}
\author{Weinan Wang}
\address{Department of Mathematics, University of Arizona, Tucson, AZ, 85712, USA}
\email{weinanwang@math.arizona.edu}

\usepackage[margin=1in]{geometry}
\usepackage{amsmath, amsthm, amssymb}
\usepackage{times}
\usepackage{color}
\usepackage{hyperref}
\newcommand{\pa}{\partial}
\newcommand{\la}{\label}
\newcommand{\fr}{\frac}
\newcommand{\na}{\nabla}
\newcommand{\be}{\begin{equation}}
\newcommand{\ee}{\end{equation}}
\newcommand{\ba}{\begin{array}{l}}
	\newcommand{\ea}{\end{array}}

\newcommand{\beg}{\begin}

\renewcommand{\l}{\Lambda}

\usepackage{MnSymbol,wasysym}

\newcommand{\R}{\mathbb R}

\def\RR{{\mathbb R}}

\date{\today}
\begin{document}
	\begin{abstract} We consider the inhomogeneous heat equation on the half-space $\mathbb R_{+}^{d}$ with Neumann boundary conditions. We prove a space-time Gevrey regularity of the solution, with a radius of analyticity uniform up to the boundary of the half-space. We also address the case of homogeneous Robin boundary conditions. Our results generalize the case of homogeneous Dirichlet boundary conditions established by Kukavica and Vicol in \cite{kv}. 
	\end{abstract}

	\maketitle

\section{Introduction}\la{intro}

We consider the heat equation
\be \la{heat}
\pa_t q - \Delta q = f 
\ee on the upper-half space
\be 
\Omega = \RR^d_+ = \left\{(x = (x_1, ..., x_{d}) \in \RR^d : x_d > 0\right\}, \quad d \geq 2
\ee with the initial condition
\be 
q(x,0) = q_0(x),
\ee and the homogeneous Neumann boundary conditions
\be \la{heatb}
\na q|_{\pa \Omega} \cdot n = 0.
\ee Here $n$ is the outward unit normal vector of $\pa \Omega$. Since $n = (0,...,0,-1)$, the condition \eqref{heatb} reduces to 
\be 
\pa_d q = 0
\ee on $\pa \Omega$, where $\pa_d$ stands for the normal derivative of $q$. The forcing term $f$ in \eqref{heat} is a function of both space and time.

Several approaches were developed over the years to study the analyticity of nonlinear parabolic equations on domains with boundaries (\cite{K, K1}) based on successive applications of the $L^2$ norms of derivatives, and without boundaries based on Fourier series techniques (\cite{ ALW, FTS, FT, LG} and references therein), a mild formulation of the complexified problem (\cite{AI}, \cite{GK}), etc. Recently, Kukavica and Vicol established in \cite{kv} a derivative reduction proof, based on classical energy inequalities, to study the analyticity up to the boundary  of the $d$-dimensional inhomogeneous  heat equation on the half-space with homogeneous Dirichlet Boundary conditions. 

In this paper, we seek a simple energy-type argument to prove the instantaneous space-time analyticity of solutions to the initial boundary value problem \eqref{heat}--\eqref{heatb}. The Neumann type of boundary conditions imposed on the solution $q$ is a source of technical difficulty, breaking down the derivative reduction approach of \cite{kv}.  More precisely, the argument in \cite{kv} uses the elliptic regularity estimate $\|u\|_{H^2} \le C\|g\|_{L^2}$ that holds for any $u$ solving the Poisson equation $\Delta u = g$ on the upper half-space with vanishing boundary conditions. In the case of Neumann boundary conditions, the $H^2$ regularity of solutions does not have that simplified form but is rather described by the bound
\be 
\|u\|_{H^2(\Omega)} \le C\left(\|g\|_{L^2(\Omega)} + \|\bar{\partial} u\|_{H^1(\Omega)} \right).
\ee This latter dependency on the $H^1$ regularity of the tangential derivatives is a constraint on deriving derivative reduction estimates analogous to \cite{kv}. We present a proof that relies on the structure of the heat equation \eqref{heat} and uses tangential interpolation inequalities rather than elliptic estimates.

For that purpose, we consider a regularity exponent $r \ge 2$, a time $T > 0$, and strictly positive small quantities $0 < \tilde{\epsilon}, \bar{\epsilon}, \epsilon \le 1$, and we define the Gevrey type norm 
\be
\beg{aligned} \la{mains}
\psi (q) &= \sum\limits_{i + j + k \ge r} \sum\limits_{|\alpha| = k} \fr{(i+j+k)^r}{(i+j+k)!} \epsilon^i \tilde{\epsilon}^j \bar{\epsilon}^{k} \|t^{i+j+k-r}\pa_t^i \pa_d^j \bar{\pa}^{\alpha} q\|_{L_{t,x}^2([0,T] \times \Omega)}
\\&\quad\quad+ \sum\limits_{i + j + k < r} \sum\limits_{|\alpha| = k} \|\pa_t^i \pa_d^j \bar{\pa}^{\alpha} q\|_{L_{t,x}^2([0,T] \times \Omega)}
\end{aligned} 
\ee where $|\alpha|$ is the sum of the components of the $(d-1)$-dimensional vector $\alpha = (\alpha_1, \dots, \alpha_{d-1})$, and $\bar{\pa}^{\alpha}$ stands for the tangential derivative $\pa_1^{\alpha_1} \pa_2^{\alpha_2} \dots \pa_{d-1}^{\alpha_{d-1}}$. All indices over which the sums are taken are assumed to be nonnegative integers.

The second sum on the right-hand side of \eqref{mains} is the $H^{r-1}([0,T] \times \Omega))$ Sobolev norm of the solution $q$ to \eqref{heat}--\eqref{heatb}, which itself is controlled by the sum of the $H^{2(r-1)}(\Omega)$ Sobolev norm of the initial data $q_0$ and the $H^{2(r-2)}([0,T] \times \Omega)$ Sobolev norm of the forcing term $f$, provided that $q_0$ obeys the compatibility conditions (see e.g. \cite{E}). 

We seek good control of the infinite sum in \eqref{mains} via a modification of the derivative reduction approach of \cite{kv}. The following theorem states our main result:

\beg{Thm} \la{T2} Let $T > 0$ and $r \ge 2$. Then there exists $\epsilon, \tilde{\epsilon} \in (0,1]$, which depend only on $T$, $r$, and $d$,  such that for any $q_0 \in H^{2r}(\Omega)$ satisfying the compatibility conditions, and $f$ sufficiently smooth, the solution $q$ of \eqref{heat}--\eqref{heatb} satisfies the estimate
\be 
\beg{aligned}
\psi &\lesssim \|q_0\|_{H^{2r}(\Omega)}  + \|f\|_{H^{2r-2}([0,T] \times \Omega)}
+ \sum\limits_{i + k \ge r-2} \sum\limits_{|\alpha| = k} \fr{(i+k+2)^r \epsilon^i \tilde{\epsilon}^{k+2} }{(i+k+2)!} \|t^{i+k+2-r} \pa_t^i \bar{\pa}^{\alpha} f\|_{L_{t,x}^2([0,T] \times \Omega)}
\\&+\quad \sum\limits_{i \ge r} \fr{(i+1)^r \epsilon^{i+1} }{(i+1)!} \|t^{i+1-r} \pa_t^i f\|_{L_{t,x}^2([0,T] \times \Omega)}
+ \sum\limits_{i \ge r-1} \fr{(i+1)^{r-1} \epsilon^{i+1}}{(i+1)!} \|t^{i+2-r}\pa_t^i \pa_df\|_{L_{t,x}^2([0,T] \times \Omega)}
\\&\quad\quad+\sum\limits_{i + j + k \ge 1+(r-3)_+} \sum\limits_{|\alpha| = k} \fr{(i+j+k+1)^{r-1} \epsilon^i \tilde{\epsilon}^{j+k+1}}{(i+j+k+1)!} \|t^{i+j+k + 2-r} \pa_t^i \pa_d^{j} \bar{\pa}^{\alpha} f  \|_{L_{t,x}^2([0,T] \times \Omega)}.
\end{aligned}
\ee  Here the notation $A \lesssim B$ means that $A \le C_{r,d}B$ for some positive constant $C_{r,d}$ depending only on $r$, the dimension $d$, and some universal constants. 
\end{Thm}

The idea of the proof is based on a decomposition of the norm \eqref{mains} into two main sums, one involving normal derivatives and one depending only on tangential and time derivatives. The terms with normal derivatives are controlled via the reduction technique of \cite{kv} in view of the fact that $\pa_dq$ solves an inhomogeneous heat equation with homogeneous Dirichlet boundary conditions. As for the sum which does not depend on the normal derivatives of solutions, we decompose it into three sub-sums, $S_1$, $S_2$, and $S_3$, where $S_1$ includes all terms with at least two tangential derivatives, $S_2$ depends on exactly one tangential derivative, and $S_3$ is the sum of the remaining time derivative terms. The estimation of $S_1$ uses the structure of the diffusion driven by $\Delta q$, which, by making use of the heat equation \eqref{heat}, allows us to reduce the number of tangential derivatives by increasing the number of normal derivatives. As for the sum $S_2$, we interpolate in the tangential variable to have an additional tangential derivative and hence have good control of $S_2$ by $S_1$. Finally, we estimate $S_3$ by reducing the number of time derivatives based on standard energy equalities. This latter reduction is mainly obtained via integration by parts, which holds even under the Neumann boundary conditions imposed on the solution.

The analogous result obtained in \cite{kv} for homogeneous Dirichlet boundary conditions was applied in \cite{CKV} to prove the Gevrey regularity of the Navier-Stokes equations on half-spaces. We believe that our result will also be useful to study the space-time analyticity of $d$-dimensional nonlinear parabolic equations with homogeneous Neumann boundary conditions (where $f = f(t,x,q)$ depends on $q$ in a nonlinear fashion).

We prove Theorem \ref{T2} in Section \ref{pt2}, and we briefly address the cases of inhomogeneous Neumann and homogeneous Robin boundary conditions in Section \ref{cr}. Throughout the paper, the letter $C$ denotes a positive universal constant that may change from line to line along the proofs. For fixed nonnegative indices $i,j,k$, we use the notation 
\be 
 \|\pa_t^i \pa_d^j \bar{\pa}^k q\|_{L_{t,x}^2([0,T] \times \Omega)}
 := \sum\limits_{|\alpha| = k} \|\pa_t^i \pa_d^j \bar{\pa}^{\alpha} q\|_{L_{t,x}^2([0,T] \times \Omega)}
\ee
for the sake of simplicity.

\section{Proof of Theorem \ref{T2}} \la{pt2}

Recalling the Gevrey norm \eqref{mains}, we decompose $\psi(q) = \psi_1(q) + \psi_2(q)$, where $\psi_1(q)$ and $\psi_2(q)$ are the following sums
\be
\beg{aligned} \la{sum1}
\psi_1(q) &=  \sum\limits_{i + j + k \ge r, j \ne 0} \fr{(i+j+k)^r}{(i+j+k)!} \epsilon^i \tilde{\epsilon}^j \bar{\epsilon}^k \|t^{i+j+k-r}\pa_t^i \pa_d^j \bar{\pa}^k q\|_{L_{t,x}^2([0,T] \times \Omega)}
\\&\quad\quad+ \sum\limits_{i + j + k < r, j \ne 0} \|\pa_t^i \pa_d^j \bar{\pa}^k q\|_{L_{t,x}^2([0,T] \times \Omega)},
\end{aligned} 
\ee and 
\be  \la{sum2}
\psi_2(q) = \sum\limits_{i + k \ge r} \fr{(i+k)^r}{(i+k)!} \epsilon^i  \bar{\epsilon}^k \|t^{i+k-r}\pa_t^i  \bar{\pa}^k q\|_{L_{t,x}^2([0,T] \times \Omega)}
+ \sum\limits_{i +  k < r} \|\pa_t^i  \bar{\pa}^k q\|_{L_{t,x}^2([0,T] \times \Omega)}.
\ee

\noindent \textbf{Estimation of $\psi_1$.} The Gevrey norm \eqref{sum1} is controlled via use of normal, tangential, and time derivative reductions. Indeed, $\psi_1$ can be rewritten as 
\be \la{kvres}
\beg{aligned}
\psi_1(q) &= \sum\limits_{i + \tilde{j} + k \ge r-1} \fr{(i + \tilde{j} + k + 1)^r}{(i+ \tilde{j} + k + 1)!} \epsilon^i \tilde{\epsilon}^{\tilde{j} + 1} \bar{\epsilon}^k \|t^{i+\tilde{j} + k - (r-1)} \pa_t^i \pa_d^{\tilde{j}} \bar{\pa}^k (\pa_dq)\|_{L^2_{t,x} ([0,T] \times \Omega)} 
\\&\quad\quad+ \sum\limits_{i + \tilde{j} + k < r-1} \|\pa_t^i \pa_d^{\tilde{j}} \bar{\pa}^k (\pa_dq)\|_{L^2_{t,x}([0,T] \times \Omega)}
\\&\le 2^{r-1} \tilde{\epsilon} \sum\limits_{i + \tilde{j} + k \ge r-1} \fr{(i + \tilde{j} + k)^{r-1}}{(i+ \tilde{j} + k )!} \epsilon^i \tilde{\epsilon}^{\tilde{j}} \bar{\epsilon}^k \|t^{i+\tilde{j} + k - (r-1)} \pa_t^i \pa_d^{\tilde{j}} \bar{\pa}^k (\pa_dq)\|_{L^2_{t,x} ([0,T] \times \Omega)} 
\\&\quad\quad+ \sum\limits_{i + \tilde{j} + k < r-1} \|\pa_t^i \pa_d^{\tilde{j}} \bar{\pa}^k (\pa_dq)\|_{L^2_{t,x}([0,T] \times \Omega)} 
\\&\le 2^{r-1} \phi(q),
\end{aligned} 
\ee where
\be 
\beg{aligned}
\phi(q) &=  \sum\limits_{i + \tilde{j} + k \ge r-1} \fr{(i + \tilde{j} + k)^{r-1}}{(i+ \tilde{j} + k )!} \epsilon^i \tilde{\epsilon}^{\tilde{j}} \bar{\epsilon}^k \|t^{i+\tilde{j} + k - (r-1)} \pa_t^i \pa_d^{\tilde{j}} \bar{\pa}^k (\pa_dq)\|_{L^2_{t,x} ([0,T] \times \Omega)} 
\\&\quad\quad+ \sum\limits_{i + \tilde{j} + k < r-1} \|\pa_t^i \pa_d^{\tilde{j}} \bar{\pa}^k (\pa_dq)\|_{L^2_{t,x}([0,T] \times \Omega)}.
\end{aligned}
\ee 
Here the first equality is obtained via the change of variable $\tilde{j} = j-1$ and the first inequality follows from an application of the algebraic inequality
\be 
\fr{(i+\tilde{j} + k + 1)^{r}}{(i+\tilde{j} + k + 1)!} = \fr{(i + \tilde{j} + k + 1)^{r-1}}{(i + \tilde{j} + k)!} \le \fr{(2(i+ \tilde{j} + k))^{r-1}}{(i + \tilde{j} + k)!} = \fr{2^{r-1} (i + \tilde{j} + k)^{r-1}}{(i + \tilde{j} + k)!} 
\ee that holds for any nonnegative integers $i, \tilde{j}, k$ whose sum is greater than or equal to 1. Also, the last inequality in \eqref{kvres} uses the boundedness of $\tilde{\epsilon}$ from above by $1$ and $r$ from below by 2. Since $\pa_dq$ solves the heat equation
\be 
\pa_t (\pa_dq) - \Delta(\pa_dq) = \pa_df
\ee with homogeneous Dirichlet boundary conditions, we can apply the derivative reduction technique of \cite{kv} and infer that there exists a positive universal constant $C$ depending on the dimension $d$ and $r$, such that for any  $C_0 \in (0, C)$, if $\epsilon$ obeys
\be \la{sum33}
\epsilon \le C_0,
\ee $\bar{\epsilon} = \bar{\epsilon}(\epsilon, T, C_0)$ obeys 
\be \la{sum38}
\fr{T \bar{\epsilon}^2}{\epsilon} + \fr{T^{\fr{1}{2}} \bar{\epsilon}}{\sqrt{\epsilon}} + T\bar{\epsilon} \le C_0,
\ee $\tilde{\epsilon} = \tilde{\epsilon} (\bar{\epsilon}, \epsilon, T, C_0)$ obeys
\be 
\fr{T \tilde{\epsilon}^2}{\epsilon} + \fr{\tilde{\epsilon}}{\epsilon} + \fr{\tilde{\epsilon}^2}{\epsilon^2} + T^2 \tilde{\epsilon}^2 + \fr{T \bar{\epsilon} \tilde{\epsilon}}{\epsilon} + \fr{T^{\fr{1}{2}}\tilde{\epsilon}}{\epsilon^{\fr{1}{2}}} + T \tilde{\epsilon} \le C_0,
\ee and 
\be \la{sum37}
0 < \tilde{\epsilon} \le \bar{\epsilon} \le \epsilon \le 1,
\ee then 
\be \la{sum3}
\phi(q) \lesssim \phi_0 + \phi_1(f)
\ee where $\phi_0$ and $\phi_1(f)$ are given by 
\be \la{sum5}
\phi_0 = \|q_0\|_{H^{2r-1}(\Omega)} + \|f\|_{H^{2r-3}([0,T] \times \Omega)}
\ee and 
\be \la{sum4}
\beg{aligned}
\phi_1(f) &= \sum\limits_{i + j + k \ge 1+(r-3)_+} \fr{(i+j+k+1)^{r-1} \epsilon^i \tilde{\epsilon}^{j+1} \bar{\epsilon}^k}{(i+j+k+1)!} \|t^{i+j+k + 2-r} \pa_t^i \pa_d^{j} \bar{\pa}^k f  \|_{L_{t,x}^2([0,T] \times \Omega)} 
\\&\quad\quad+ \sum\limits_{i+k \ge (r-3)_+} \fr{(i+k+2)^{r-1} \epsilon^i \bar{\epsilon}^{k+2}}{(i+k+2)!} \|t^{i+k+3 - r} \pa_t^i \pa_d\bar{\pa}^k f\|_{L_{t,x}^2([0,T] \times \Omega)}
\\&\quad\quad\quad\quad+ \sum\limits_{i \ge r-1} \fr{(i+1)^{r-1} \epsilon^{i+1}}{(i+1)!} \|t^{i+2-r}\pa_t^i \pa_df\|_{L_{t,x}^2([0,T] \times \Omega)},
\end{aligned}
\ee provided that the quantities \eqref{sum5} and \eqref{sum4} are finite. We refer the reader to \cite[(4.9)--(4.12)]{kv} for the choices of $\epsilon, \tilde{\epsilon}, \bar{\epsilon}$ given by \eqref{sum33}--\eqref{sum37}. Further assumptions on $\epsilon, \tilde{\epsilon}, \bar{\epsilon}$ will be imposed later.  

\noindent \textbf{Estimation of $\psi_2$.} Here $\psi_2$ is the trickier term. We split $\psi_2$ into $\psi_2 = \psi_{2,1} + \psi_{2,2}$ where 
\be \la{sum7}
\psi_{2,1} = \sum\limits_{i +  k < r} \|\pa_t^i  \bar{\pa}^k q\|_{L_{t,x}^2([0,T] \times \Omega)}
\ee and 
\be \la{sum8}
\psi_{2,2} = \sum\limits_{i + k \ge r} \fr{(i+k)^r}{(i+k)!} \epsilon^i  \bar{\epsilon}^k \|t^{i+k-r}\pa_t^i  \bar{\pa}^k q\|_{L_{t,x}^2([0,T] \times \Omega)}.
\ee The norm $\psi_{2,1}$ is bounded by the sum of the $H^{2(r-1)}(\Omega)$ norm of $q_0$ and the $H^{2(r-2)}([0,T] \times \Omega))$ norm of $f$. In order to control $\psi_{2,2}$, we perform tangential and time derivative reductions. However, we do not appeal to elliptic estimates but use the PDE obeyed by $q$ instead. We decompose $\psi_{2,2}$ into the sum $\psi_{2,2} = \psi_{2,2,1} + \psi_{2,2,2} + \psi_{2,2,3}$, where
\be \la{MMM}
\psi_{2,2,1} = \sum\limits_{i+k \ge r, k \ge 2} \fr{(i+k)^r}{(i+k)!} \epsilon^i  \bar{\epsilon}^k \|t^{i+k-r}\pa_t ^i \bar{\pa}^k q\|_{L_{t,x}^2([0,T] \times \Omega)},
\ee
\be 
\psi_{2,2,2} = \sum\limits_{i \ge r-1} \fr{(i+1)^r}{(i+1)!} \epsilon^i  \bar{\epsilon} \|t^{i+1-r}\pa_t^i  \bar{\pa} q\|_{L_{t,x}^2([0,T] \times \Omega)},
\ee and
\be 
\psi_{2,2,3} = \sum\limits_{i \ge r}  \fr{i^r}{i!} \epsilon^i  \|t^{i-r}\pa_t^i  q\|_{L_{t,x}^2([0,T] \times \Omega)}.
\ee
We start by estimating $\psi_{2,2,1}$. The main idea is to reduce the number of horizontal derivatives by increasing the number of vertical derivatives, which, eventually, allows us to control $\psi_{2,2,1}$ by the sum $\psi_1$. A loss of two tangential derivatives is equivalent to a gain of two normal derivatives, a fact that is based on the heat equation \eqref{heat}, from which we obtain the relation
\be \la{neuf}
\Delta_{d-1} q = - \pa_d\pa_dq + \pa_t q  - f.
\ee Here $\Delta_{d-1}$ stands for the $(d-1)$-dimensional Laplace operator, 
\be 
\Delta_{d-1} q := \pa_1\pa_1 q + \dots + \pa_{d-1}\pa_{d-1} q.
\ee 
The norm $\|\cdot\|_{L_{t,x}^2([0,T] \times \Omega)}$ is equivalent to
\be 
\|\cdot\|_{L_{t,x}^2([0,T] \times \Omega)} = \left\| \|\cdot\|_{L_{x_1, \dots, x_{d-1}}^2(\R^{d-1})} \right\|_{L_{t, x_d}^2 ([0,T] \times (0, \infty))},
\ee and so the sum $\psi_{2,2,1}$ can be written as 
\be 
\psi_{2,2,1} = \sum\limits_{i+k \ge r, k \ge 2} \fr{(i+k)^r}{(i+k)!} \epsilon^i \bar{\epsilon}^k \left\|t^{i+k-r} \|\pa_t^i \bar{\pa}^k q\|_{L_{x_1, \dots, x_{d-1}}^2(\R^{d-1})} \right\|_{L_{t, x_d}^2 ([0,T] \times (0, \infty))}.
\ee  Denoting the inverse of the square root of the Laplace operator $-\Delta_{d-1}$ on the whole space $\R^{d-1}$ by $\l_{d-1}^{-1}$, and exploiting the boundedness of the Riesz transform operator $\na_{d-1} \l_{d-1}^{-1}$ on $L^2(\RR^{d-1})$, we have the following $(d-1)$-dimensional elliptic regularity estimate 
\be \la{riesz}
\|\pa_{x_s} \pa_{x_r} \rho\|_{L^2(\RR^{d-1})} = \|\pa_{x_s} \l_{d-1}^{-1} \pa_{x_r} \l_{d-1}^{-1} \Delta_{d-1} \rho\|_{L^2(\RR^{d-1})} \le C\|\Delta_{d-1} \rho\|_{L^2(\RR^{d-1})}
\ee for any $\rho \in H^2(\RR^{d-1})$, and any $s, r \in \left\{1, \dots, d-1\right\}$. We point out that the first equality in \eqref{riesz} follows from the fact the operators $\na_{d-1}$ and $\l_{d-1}^{-1}$ are Fourier multipliers in the whole space setting, so they commute.  Accordingly, for any nonnegative integers $i \ge 0$ and $k \ge 2$, we have
\be 
\|\pa_t^i \bar{\pa}^k q\|_{L_{x_1, \dots, x_{d-1}}^2(\R^{d-1})}
\le \|\pa_t^i \bar{\pa}^{k-2} q\|_{\dot{H}^2(\RR^{d-1})}
\le C\|\pa_t^i \bar{\pa}^{k-2} \Delta_{d-1} q\|_{L_{x_1, \dots, x_{d-1}}^2 (\mathbb R^{d-1})}
\ee which, followed by an application of the relation \eqref{neuf}, yields the estimate
\be \la{sum6}
\beg{aligned} 
\psi_{2,2,1}
&\le  \sum\limits_{i+k \ge r, k \ge 2} \fr{(i+k)^r}{(i+k)!} \epsilon^i  \bar{\epsilon}^k \|t^{i+k-r}\pa_t ^i \pa_d^2 \bar{\pa}^{k-2}  q\|_{L_{t,x}^2([0,T] \times \Omega)}
\\&\quad\quad+  \sum\limits_{i+k \ge r, k \ge 2} \fr{(i+k)^r}{(i+k)!} \epsilon^i  \bar{\epsilon}^k \|t^{i+k-r}\pa_t ^{i+1} \bar{\pa}^{k-2}  q\|_{L_{t,x}^2([0,T] \times \Omega)}
\\&\quad\quad\quad\quad+  \sum\limits_{i+k \ge r -2} \fr{(i+k+2)^r}{(i+k+2)!} \epsilon^i  \bar{\epsilon}^{k+2} \|t^{i+k+2-r}\pa_t ^i \bar{\pa}^{k} f\|_{L_{t,x}^2([0,T] \times \Omega)}.
\end{aligned}
\ee
The first sum in \eqref{sum6} is controlled by a constant multiple of the Gevrey norm $\psi_1(q)$ given by \eqref{sum1} due to the presence of second-order normal derivatives. In other words, we have 
\be 
\beg{aligned}
&\sum\limits_{i+k \ge r, k \ge 2} \fr{(i+k)^r}{(i+k)!} \epsilon^i  \bar{\epsilon}^k \|t^{i+k-r}\pa_t ^i \pa_d^2 \bar{\pa}^{k-2}  q\|_{L_{t,x}^2([0,T] \times \Omega)}
\\&= \fr{\bar{\epsilon}^2}{\tilde{\epsilon}^2} \sum\limits_{i+k +2 \ge r} \fr{(i+2 + k)^r}{(i+2 + k)!} \epsilon^i \tilde{\epsilon}^2 \bar{\epsilon}^k \|t^{i+2 + k-r}\pa_t ^i \pa_d^2 \bar{\pa}^{k}  q\|_{L_{t,x}^2([0,T] \times \Omega)}
\\&\le \fr{\bar{\epsilon}^2}{\tilde{\epsilon}^2} \psi_1(q)
\le \fr{\bar{\epsilon}^2}{\tilde{\epsilon}^2} 2^{r-1} \left(\phi_0 + \phi_1(f)\right) = 2^{r-1} \left(\phi_0 + \phi_1(f)\right),
\end{aligned}
\ee provided that $\bar{\epsilon} = \tilde{\epsilon}$. We recall that $\phi_0$ and $\phi_1(f)$ are given by \eqref{sum5} and \eqref{sum4} respectively. 
The second sum in \eqref{sum6} is bounded by a small constant multiple of \eqref{sum8} up to an additive constant depending only on $q_0$ and $f$. Indeed, 
\be 
\beg{aligned}
&\sum\limits_{i+k \ge r, k \ge 2} \fr{(i+k)^r}{(i+k)!} \epsilon^i  \bar{\epsilon}^k \|t^{i+k-r}\pa_t ^{i+1} \bar{\pa}^{k-2}  q\|_{L_{t,x}^2([0,T] \times \Omega)}
\\&= \fr{\bar{\epsilon}^2}{\epsilon} \sum\limits_{i+k \ge r-1} \fr{(i+k+1)^r}{(i+k+1)!} \epsilon^i  \bar{\epsilon}^k \|t^{i+k-r+ 1}\pa_t ^{i} \bar{\pa}^{k}  q\|_{L_{t,x}^2([0,T] \times \Omega)}
\\&\le \fr{T\bar{\epsilon}^2}{\epsilon}  \sum\limits_{i+k \ge r} \fr{(i+k+1)^{r-1}}{(i+k)!} \epsilon^i  \bar{\epsilon}^k \|t^{i+k-r}\pa_t ^{i} \bar{\pa}^{k}  q\|_{L_{t,x}^2([0,T] \times \Omega)}
+ \fr{\bar{\epsilon}^2}{\epsilon}\fr{r^r}{r!} \sum\limits_{i+k = r-1}  \epsilon^i  \bar{\epsilon}^k \|\pa_t ^{i} \bar{\pa}^{k}  q\|_{L_{t,x}^2([0,T] \times \Omega)}
\\&\le \fr{2^{r-1} T\bar{\epsilon}^2}{\epsilon}  \sum\limits_{i+k \ge r} \fr{(i+k)^{r}}{(i+k)!} \epsilon^i  \bar{\epsilon}^k \|t^{i+k-r}\pa_t ^{i} \bar{\pa}^{k}  q\|_{L_{t,x}^2([0,T] \times \Omega)}
+ \fr{r^r}{r!}  \sum\limits_{i+k = r-1}  \|\pa_t ^{i} \bar{\pa}^{k}  q\|_{L_{t,x}^2([0,T] \times \Omega)}
\end{aligned}
\ee by \eqref{sum37}. In view of the condition \eqref{sum38}, we infer that
\be  \la{sum66}
\beg{aligned}
&\sum\limits_{i+k \ge r, k \ge 2} \fr{(i+k)^r}{(i+k)!} \epsilon^i  \bar{\epsilon}^k \|t^{i+k-r}\pa_t ^{i+1} \bar{\pa}^{k-2}  q\|_{L_{t,x}^2([0,T] \times \Omega)}
\\&\le 2^{r-1} C_0 \psi_{2,2} + C\|q_0\|_{H^{2(r-1)(\Omega)}} + C\|f\|_{H^{2(r-2)}([0,T] \times \Omega)}
\\&\le \delta \psi_{2,2} + C\|q_0\|_{H^{2(r-1)(\Omega)}} + C\|f\|_{H^{2(r-2)}([0,T] \times \Omega)}
\end{aligned}
\ee  provided that $C_0$ is chosen to be smaller than $\fr{\delta}{2^{r-1}}$. Here $\delta$ is a positive constant that will be determined later. 

Now we proceed to estimate $\psi_{2,2,2}$. The main idea is to increase the number of tangential derivative by one via interpolation and show that $\psi_{2,2,2}$ is dominated by the sum of $\psi_{2,2,1}$ and $\psi_{2,2,3}$, reducing consequently the problem to a time derivative reduction. We need the following elementary lemma:

\beg{lem} \la{elem} Let $i \ge 1$ and $r \ge 1$ be some integers. Then the following estimate 
\be 
\fr{(i+1)^r}{(i+1)!} \le \sqrt{2^{r+1}} \sqrt{\fr{(i+2)^r}{(i+2)!}} \sqrt{\fr{i^r}{i!}}
\ee holds.
\end{lem}

\textbf{Proof of Lemma \ref{elem}.} We have  
\be 
\beg{aligned}
&\sqrt{\fr{i! (i+2)!}{[(i+1)!]^2}} \sqrt{\fr{(i+1)^{2r}}{i^r (i+2)^r}}
= \sqrt{\fr{i+2}{i+1}} \sqrt{\left(\fr{i+1}{i} \right)^r \left(\fr{i+1}{i+2} \right)^r}
\\&= \sqrt{1 + \fr{1}{i+1}} \sqrt{\left( 1+ \fr{1}{i}\right)^r \left(\fr{i+1}{i+2} \right)^r}
\le \sqrt{2} \sqrt{2^r} = \sqrt{2^{r+1}}.
\end{aligned}
\ee

In view of Lemma \ref{elem} and the $(d-1)$-dimensional interpolation inequality 
\be 
\|\bar{\pa} \rho\|_{L_{x_1, \dots, x_{d-1}}^2(\R^{d-1})} \lesssim \|\bar{\pa} \bar{\pa} \rho\|_{L_{x_1, \dots, x_{d-1}}^2(\R^{d-1})}^{\fr{1}{2}} \|\rho\|_{L_{x_1, \dots, x_{d-1}}^2(\R^{d-1})}^{\fr{1}{2}} + \|\rho\|_{L_{x_1, \dots, x_{d-1}}^2(\R^{d-1})}
\ee that holds for any $\rho \in H^2(\R^{d-1})$, we have 
\be 
\beg{aligned}
\psi_{2,2,2} 
&= \sum\limits_{i \ge r-1} \fr{(i+1)^r}{(i+1)!}\epsilon^i \bar{\epsilon} \|t^{i+1-r} \pa_t^i \bar{\pa}q\|_{L_{t,x}^2 ([0,T] \times \Omega)}
\\&=\fr{r^r}{r!} \epsilon^{r-1} \bar{\epsilon} \|\pa_t^{r-1} \bar{\pa} q\|_{L_{t,x}^2 ([0,T] \times \Omega)}
+ \sum\limits_{i \ge r} \fr{(i+1)^r}{(i+1)!}\epsilon^i \bar{\epsilon} \|t^{i+1-r} \pa_t^i \bar{\pa}q\|_{L_{t,x}^2 ([0,T] \times \Omega)}
\\&\le C \left(\|q_0\|_{H^{2r-1}(\Omega)} + \|f\|_{H^{2r-3}([0,T] \times \Omega)} \right) 
+ \sum\limits_{i \ge r} \fr{(i+1)^r}{(i+1)!} \epsilon^i \bar{\epsilon} \|t^{i+1-r} \pa_t^i q\|_{L_{t,x}^2([0,T] \times \Omega)}
\\&\quad+ C\sqrt{2^{r+1}} \sum\limits_{i \ge r} \sqrt{\fr{(i+2)^r}{(i+2)!}} \sqrt{\fr{i^r}{i!}} \epsilon^i \bar{\epsilon} \left\| \|\bar{\pa}\bar{\pa} (t^{i+1-r} \pa_t^i q) \|_{L^2(\R^{d-1})}^{\fr{1}{2}} \|t^{i+1-r} \pa_t^i q \|_{L^2(\R^{d-1})}^{\fr{1}{2}} \right\|_{L_{t, x_d}^2}
\\&= C \left(\|q_0\|_{H^{2r-1}(\Omega)} + \|f\|_{H^{2r-3}([0,T] \times \Omega)} \right) 
+ \sum\limits_{i \ge r} \fr{(i+1)^r}{(i+1)!} \epsilon^i \bar{\epsilon} \|t^{i+1-r} \pa_t^i q\|_{L_{t,x}^2([0,T] \times \Omega)}
\\&\quad+ C\sqrt{2^{r+1}} \sum\limits_{i \ge r}   \left\| \left\| \fr{(i+2)^r}{(i+2)!} \epsilon^i \bar{\epsilon}^2  \bar{\pa}\bar{\pa} (t^{i+2-r} \pa_t^i q) \right\|_{L^2(\R^{d-1})}^{\fr{1}{2}} \left\| \fr{i^r}{i!} \epsilon^i t^{i-r}  \pa_t^i q \right\|_{L^2(\R^{d-1})}^{\fr{1}{2}} \right\|_{L_{t, x_d}^2},
\end{aligned}
\ee which, after using Young's inequality, boils down to
\be 
\beg{aligned}
&\psi_{2,2,2}
\le C \left(\|q_0\|_{H^{2r-1}(\Omega)} + \|f\|_{H^{2r-3}([0,T] \times \Omega)} \right) 
+ \sum\limits_{i \ge r} \fr{(i+1)^r}{(i+1)!} \epsilon^i \bar{\epsilon} \|t^{i+1-r} \pa_t^i q\|_{L_{t,x}^2([0,T] \times \Omega)}
\\&+ C\sqrt{2^{r-1}} \sum\limits_{i \ge r}   \fr{(i+2)^r}{(i+2)!} \epsilon^i \bar{\epsilon}^2 \|t^{i+2-r} \pa_t^i \bar{\pa}\bar{\pa}  q\|_{L_{t,x}^2([0,T] \times \Omega)}
+ C\sqrt{2^{r-1}} \sum\limits_{i \ge r}   \fr{i^r}{i!}  \epsilon^i \|t^{i-r}  \pa_t^i q\|_{L_{t,x}^2([0,T] \times \Omega)}
\\&\le C\|q_0\|_{H^{2r-1}(\Omega)} + C\|f\|_{H^{2r-3}([0,T] \times \Omega)} 
+ C2^{r-1} \bar{\epsilon} T  \psi_{2,2,3} + C\sqrt{2^{r-1}}\psi_{2,2,1} + C\sqrt{2^{r-1}} \psi_{2,2,3}.
\end{aligned}
\ee Since $T\bar{\epsilon} \le C_0$, we obtain
\be 
\beg{aligned} \la{sum666}
\psi_{2,2,2} &\le C\|q_0\|_{H^{2r-1}(\Omega)} + C\|f\|_{H^{2r-3}(([0,T] \times \Omega))} + C_1\sqrt{2^{r-1}}\psi_{2,2,1} + C\psi_{2,2,3},
\end{aligned} 
\ee  where $C$ is a positive constant depending only on $r$ and $d$, and $C_1$ is a positive universal constant.

We end the proof by estimating $\psi_{2,2,3}$. We take the scalar product in $L^2$ of the heat equation \eqref{heat} with $\pa_t q$. We integrate by parts the diffusion term $(-\Delta q, \pa_t q)_{L^2}$ using the homogeneous Neumann boundary conditions. We apply the Cauchy-Schwarz inequality to bound the forcing term $(f, \pa_t q)_{L^2}$ and then make use of Young's inequality to obtain the energy inequality
\be \la{ddd}
\|\pa_t q\|_{L^2}^2 + \fr{d}{dt} \|\na q\|_{L^2}^2 \le \|f\|_{L^2}^2,
\ee which yields 
\be 
\int_{0}^{T} \|\pa_t q\|_{L^2}^2 dt + \|\na q(T)\|_{L^2}^2 \le \|\na q(0)\|_{L^2} + \int_{0}^{T} \|f(t)\|_{L^2}^2 dt
\ee after integrating in time from $0$ to $T$. This latter estimate holds for any inhomogeneous heat equation with homogeneous Neumann boundary conditions. Consequently, it applies to the equation
\be \la{rr11}
\pa_t (t^{i-r} \pa_t^{i-1} q) - \Delta (t^{i-r} \pa_t^{i-1} q) = (i-r)t^{i-1-r} \pa_t^{i-1}q + t^{i-r}\pa_t^{i-1}f
\ee for any $i \ge r+1$. Since $t^{i-r}\pa_t^{i-1}q$ vanishes at the initial time $t=0$, we obtain
\be 
\|\pa_t(t^{i-r} \pa_t^{i-1}q)\|_{L_{t,x}^2([0,T] \times \Omega)} \le \|(i-r)t^{i-1-r}\pa_t^{i-1} q\|_{L_{t,x}^2([0,T] \times \Omega)}  + \|t^{i-r} \pa_t^{i-1} f\|_{L_{t,x}^2([0,T] \times \Omega)} ,
\ee which boils down to 
\be \la{timereduction}
\|t^{i-r} \pa_t^{i} q\|_{L_{t,x}^2([0,T] \times \Omega)} \lesssim \|(i-r)t^{i-1-r}\pa_t^{i-1} q\|_{L_{t,x}^2([0,T] \times \Omega)}  + \|t^{i-r} \pa_t^{i-1} f\|_{L_{t,x}^2([0,T] \times \Omega)} 
\ee for any $i \ge r+1$. By making use of \eqref{timereduction}, we estimate $\psi_{2,2,3}$ as follows,
\be 
\beg{aligned}
\psi_{2,2,3} &\lesssim \fr{r^r}{r!} \epsilon^r \|\pa_t^r q\|_{L_{t,x}^2([0,T] \times \Omega)}
+ \sum\limits_{i \ge r+1} \fr{i^r}{i!}\epsilon^i  (i-r) \|t^{i-1 - r} \pa_t^{i-1} q\|_{L_{t,x}^2([0,T] \times \Omega)}
\\&\quad\quad+ \sum\limits_{i \ge r+1} \fr{i^r}{i!} \epsilon^i \|t^{i-r} \pa_t^{i-1} f\|_{L_{t,x}^2([0,T] \times \Omega)}
\\&\le C (\|q_0\|_{H^{2r}(\Omega)} + \|f\|_{H^{2r-2}([0,T] \times \Omega)}) 
+ C_2 2^{r} \epsilon \psi_{2,2,3} 
\\&\quad\quad+ C\sum\limits_{i \ge r} \fr{(i+1)^r}{(i+1)!} \epsilon^{i+1} \|t^{i+1-r} \pa_t^{i} f\|_{L_{t,x}^2([0,T] \times \Omega)},
\end{aligned}
\ee from which we obtain 
\be \la{sum6666}
\psi_{2,2,3} \le C\left(\|q_0\|_{H^{2r}(\Omega)} + \|f\|_{H^{2r-2}([0,T] \times \Omega)}
 + \sum\limits_{i \ge r} \fr{(i+1)^r}{(i+1)!} \epsilon^{i+1} \|t^{i+1-r} \pa_t^{i} f\|_{L_{t,x}^2([0,T] \times \Omega)}\right)
\ee provided that $\epsilon \le C_0 \le \fr{1}{C_22^{r+1}}$. 
Putting \eqref{sum6}--\eqref{sum66}, \eqref{sum666} and \eqref{sum6666} together, we conclude that
\be 
\beg{aligned}
\psi_{2,2} &\le C\left(\|q_0\|_{H^{2r}(\Omega)} + \|f\|_{H^{2r-2}([0,T] \times \Omega)}\right)
+ \delta \psi_{2,2} + C_1 \sqrt{2^{r-1}} \delta \psi_{2,2} 
\\&\quad\quad+ C\phi_1(f) +  C \sum\limits_{i \ge r} \fr{(i+1)^r}{(i+1)!} \epsilon^{i+1} \|t^{i+1-r} \pa_t^{i} f\|_{L_{t,x}^2([0,T] \times \Omega)}
\\&\quad\quad\quad\quad+ C \sum\limits_{i+k \ge r -2} \fr{(i+k+2)^r}{(i+k+2)!} \epsilon^i  \bar{\epsilon}^{k+2} \|t^{i+k+2-r}\pa_t ^i \bar{\pa}^{k} f\|_{L_{t,x}^2([0,T] \times \Omega)}.
\end{aligned}
\ee for any $\delta > 0$. Choosing $\delta$ so that $\delta (1+C_1 \sqrt{2^{r-1}}) \le \fr{1}{2}$, we obtain the following bound for $\psi_{2,2}$, 
\be 
\beg{aligned}
\psi_{2,2} &\le C\left(\|q_0\|_{H^{2r}(\Omega)} + \|f\|_{H^{2r-2}([0,T] \times \Omega)}\right)
+  C \sum\limits_{i \ge r} \fr{(i+1)^r}{(i+1)!} \epsilon^{i+1} \|t^{i+1-r} \pa_t^{i} f\|_{L_{t,x}^2([0,T] \times \Omega)}
\\&\quad\quad\quad\quad+ C \sum\limits_{i+k \ge r -2} \fr{(i+k+2)^r}{(i+k+2)!} \epsilon^i  \bar{\epsilon}^{k+2} \|t^{i+k+2-r}\pa_t ^i \bar{\pa}^{k} f\|_{L_{t,x}^2([0,T] \times \Omega)} + C\phi_1(f).
\end{aligned}
\ee
Therefore, we infer that 
\be 
\beg{aligned}
\psi &\lesssim \|q_0\|_{H^{2r}(\Omega)}  + \|f\|_{H^{2r-2}([0,T] \times \Omega)}
+ \sum\limits_{i + k \ge r-2} \fr{(i+k+2)^r}{(i+k+2)!} \epsilon^i \bar{\epsilon}^{k+2} \|t^{i+k+2-r} \pa_t^i \bar{\pa}^k f\|_{L_{t,x}^2([0,T] \times \Omega)}
\\&+\quad \sum\limits_{i \ge r} \fr{(i+1)^r}{(i+1)!}\epsilon^{i+1} \|t^{i+1-r} \pa_t^i f\|_{L_{t,x}^2([0,T] \times \Omega)}
\\&\quad\quad+\sum\limits_{i + j + k \ge 1+(r-3)_+} \fr{(i+j+k+1)^{r-1} \epsilon^i \tilde{\epsilon}^{j+1} \bar{\epsilon}^k}{(i+j+k+1)!} \|t^{i+j+k + 2-r} \pa_t^i \pa_d^{j} \bar{\pa}^k f  \|_{L_{t,x}^2([0,T] \times \Omega)} 
\\&\quad\quad\quad+ \sum\limits_{i+k \ge (r-3)_+} \fr{(i+k+2)^{r-1} \epsilon^i \bar{\epsilon}^{k+2}}{(i+k+2)!} \|t^{i+k+3 - r} \pa_t^i \pa_d\bar{\pa}^k f\|_{L_{t,x}^2([0,T] \times \Omega)}
\\&\quad\quad\quad\quad+ \sum\limits_{i \ge r-1} \fr{(i+1)^{r-1} \epsilon^{i+1}}{(i+1)!} \|t^{i+2-r}\pa_t^i \pa_df\|_{L_{t,x}^2([0,T]\times \Omega)}
\end{aligned}
\ee for any $\epsilon, \tilde{\epsilon}, \bar{\epsilon} \in (0,1]$ obeying conditions \eqref{sum33}--\eqref{sum37}, together with $\tilde{\epsilon} = \bar{\epsilon}$. This finishes the proof of Theorem~\ref{T2}.

\section{Remarks on the inhomogeneous Neumann and homogeneous Robin boundary conditions} \la{cr}

Theorem \ref{T2} can be generalized to the case of general Neumann boundary conditions: 

\beg{rem} Let $g$ be a time-independent sufficiently smooth function defined on $\RR^{d-1}$. Solutions to the inhomogeneous heat equation on the half space $\Omega = \RR^d_+$ with general Neumann boundary conditions
\be \la{inneu}
\pa_d q = g
\ee are analytic in space and time, a fact that follows by adapting the approach of Theorem \ref{T2} to the boundary value problem formed by \eqref{heat} and \eqref{inneu}. Indeed, the norm \eqref{mains} can be decomposed into two sub-sums $S_1$ and $S_2$, where $S_1$ encompasses all normal derivatives and $S_2$ depends only on tangential and time derivatives. Since the function $v:= \pa_d q -g$ solves the heat equation
\be 
\pa_t (\pa_d q - g) - \Delta (\pa_d q -g) = \pa_d f + \bar{\pa} \cdot \bar{\pa} g
\ee and vanishes on the boundary of $\Omega$, then we obtain good control of $S_1$ by the Sobolev norm of the initial datum and the Gevrey norms of both $f$ and $g$. Here we abused notation and wrote $g$ for the extension $\tilde{g} (x_1, \dots, x_d) = g(x_1, \dots, x_{d-1})$. As for the sum $S_2$, we perform the same decomposition
strategy as for $\psi_2$ in the proof of Theorem \ref{T2}, implementing henceforth our idea of decreasing the number of tangential derivatives by increasing the number of normal derivatives. The one and only main difference resides in the derivation of the energy inequality \eqref{ddd}, which relies on integration by parts and use of the homogeneous type of Neumann boundary conditions. In the case of \eqref{inneu}, the following analogous ordinary differential equation holds
\be  \la{rr12}
\|\pa_t q\|_{L^2}^2 - \int_{\Omega} \Delta q \pa_t q dx = \int_{\Omega} f \pa_t q dx
\ee which, due to H\"older and Young inequalities, reduces to 
\be 
\|\pa_t q \|_{L^2}^2 +  \fr{d}{dt} \|\na q - G\|_{L^2}^2 \le \|f\|_{L^2}^2,
\ee  with $G = (0, \dots, 0, - \tilde{g})$. Here we used
\be 
\beg{aligned}
&\fr{1}{2} \fr{d}{dt} \|\na q - G\|_{L^2}^2 
= \int_{\Omega} (\na q - G) \cdot \pa_t (\na q - G) dx
\\&= \int_{\Omega} (\na q - G) \cdot \pa_t \na q dx
= -\int_{\Omega} \na \cdot (\na q - G) \cdot \pa_t q dx
= -\int_{\Omega} \Delta q  \pa_t q dx
\end{aligned}
\ee that holds in view of the divergence-free condition obeyed by $G$, and the vanishing property $(\na q - G)|_{\pa \Omega} \cdot n = 0$ for $n = (0, \dots, 0 , -1)$. However, the energy equality \eqref{rr12} is not needed to perform time derivative reduction, as we seek bounds for the solution of the heat equation \eqref{rr11} obeyed by $t^{i - r}\pa_t^{i-1} q$, which has a vanishing normal derivative on the boundary of the half-space for any $i \ge r+1$. The details follow along the lines of the proof of Theorem \ref{T2} and will be omitted.
\end{rem}

Our approach also applies to obtain the space-time analyticity of the heat equation with homogeneous Robin boundary conditions:

\beg{rem} The Gevrey regularity of solutions to the inhomogeneous heat equation \eqref{heat} on the half space $\Omega = \R^d_+$ with homogeneous Robin boundary conditions
\be 
(aq + b\pa_d q) |_{\pa \Omega} = 0
\ee reduces to a question of Sobolev global regularity, under some conditions imposed on $a$ and $b$. Indeed, if $a=0$ or $b=0$, then the problem boils down to the case of homogeneous Neumann or homogeneous Dirichlet boundary conditions. If both $a$ and $b$ are nonvanishing and have the same sign, then we repeat the same strategy of Theorem \ref{T2} and decompose the norm \eqref{mains} into two sub-sums, $S_1$ involving the vertical derivative components and $S_2$ involving only horizontal and time derivatives. As $v = \fr{a}{b} q + \pa_d q$ solves the heat equation 
\be 
\pa_t \left(\fr{a}{b} q + \pa_d q \right) - \Delta \left(\fr{a}{b} q + \pa_d q \right) = \fr{a}{b} f + \pa_d f
\ee with homogeneous Dirichlet boundary conditions, we can bound $S_1$ by 
\be \la{rr2}
\beg{aligned}
S_1&\le  2^{r-1} \tilde{\epsilon} \sum\limits_{i + j + k \ge r-1} \fr{(i + j + k)^{r-1}}{(i+ j + k )!} \epsilon^i \tilde{\epsilon}^{j} \bar{\epsilon}^k \|t^{i+j + k - (r-1)} \pa_t^i \pa_d^{j} \bar{\pa}^k (\pa_dq) \|_{L^2_{t,x} ([0,T] \times \Omega)}
\\&\quad\quad+ \sum\limits_{i + j + k < r-1} \|\pa_t^i \pa_d^{j} \bar{\pa}^k \pa_d q\|_{L^2_{t,x}([0,T] \times \Omega)} 
\\&\le 2^{r-1} \tilde{\epsilon} \sum\limits_{i + j + k \ge r-1} \fr{(i + j + k)^{r-1}}{(i+ j + k )!} \epsilon^i \tilde{\epsilon}^{j} \bar{\epsilon}^k \|t^{i+j + k - (r-1)} \pa_t^i \pa_d^{j} \bar{\pa}^k \left(\pa_dq + \fr{a}{b}q \right)\|_{L^2_{t,x} ([0,T] \times \Omega)} 
\\&\quad\quad+ 2^{r-1} \tilde{\epsilon} \left|\fr{a}{b} \right| \sum\limits_{i + j + k \ge r-1} \fr{(i + j + k)^{r-1}}{(i+ j + k )!} \epsilon^i \tilde{\epsilon}^{j} \bar{\epsilon}^k \|t^{i+j + k - (r-1)} \pa_t^i \pa_d^{j} \bar{\pa}^k q\|_{L^2_{t,x} ([0,T] \times \Omega)} 
\\&\quad\quad+ \sum\limits_{i + j + k < r-1} \|\pa_t^i \pa_d^{j} \bar{\pa}^k \pa_d q\|_{L^2_{t,x}([0,T] \times \Omega)} 
\end{aligned}
\ee as shown in \eqref{kvres}, and obtain control of the first sum in the last inequality by applying the result of \cite{kv}. Regarding the second sum in \eqref{rr2}, it can be controlled as follows, 
\be 
\beg{aligned}
&2^{r-1} \tilde{\epsilon} \left|\fr{a}{b} \right| \sum\limits_{i + j + k \ge r-1} \fr{(i + j + k)^{r-1}}{(i+ j + k )!} \epsilon^i \tilde{\epsilon}^{j} \bar{\epsilon}^k \|t^{i+j + k - (r-1)} \pa_t^i \pa_d^{j} \bar{\pa}^k q\|_{L^2_{t,x} ([0,T] \times \Omega)} 
\\&\le 2^{r-1} \tilde{\epsilon} \left|\fr{a}{b} \right| \fr{(r-1)^{r-1}}{(r-1)!} \sum\limits_{i + j + k = r-1} \epsilon^i \tilde{\epsilon}^{j} \bar{\epsilon}^k \| \pa_t^i \pa_d^{j} \bar{\pa}^k q\|_{L^2_{t,x} ([0,T] \times \Omega)} 
\\&\quad\quad+ 2^{r-1} T \tilde{\epsilon} \left|\fr{a}{b} \right| \sum\limits_{i + j + k \ge r} \fr{(i + j + k)^{r}}{(i+ j + k )!} \epsilon^i \tilde{\epsilon}^{j} \bar{\epsilon}^k \|t^{i+j + k -r} \pa_t^i \pa_d^{j} \bar{\pa}^k q\|_{L^2_{t,x} ([0,T] \times \Omega)} 
\end{aligned}
\ee where the first sum is bounded by the $H^{r-1}([0,T] \times \Omega)$ Sobolev norm of the solution, and the second sum is bounded by a small constant multiple of the Gevrey norm \eqref{mains}, provided that $2^{r-1} \tilde{\epsilon} \left|\fr{a}{b} \right|$ is sufficiently small.  The sum $S_2$ is treated as $\psi_2$ in the proof of Theorem \ref{T2}, but an ODE analogous to \eqref{ddd} is needed. In fact,  we have
\be  \la{rr1}
\|\pa_t q\|_{L^2}^2 + \int_{\Omega} \na q \cdot \pa_t \na q dx - \int_{\pa \Omega} \pa_t q \pa_d q d\sigma(x) = \int_{\Omega} f \pa_t q dx
\ee which, after use of the Robin boundary conditions, reduces to 
\be 
\|\pa_t q\|_{L^2}^2 + \int_{\Omega} \na q \cdot \pa_t \na q dx + \fr{a}{b} \int_{\pa \Omega} q \pa_t q  d\sigma(x) = \int_{\Omega} f \pa_t q dx,
\ee and so
\be 
\|\pa_t q\|_{L^2}^2 + \fr{1}{2} \fr{d}{dt} \left(\|\na q\|_{L^2}^2 + \fr{a}{b} \int_{\pa \Omega} q^2 d\sigma(x) \right) = \int_{\Omega} f \pa_t q dx.
\ee Now we proceed as for the proof of Theorem \ref{T2}. We omit further details. 
\end{rem}

\section{Acknowledgment}

WW was partially supported by an AMS-Simons travel grant.

 \bibliographystyle{abbrv}
 \bibliography{AW2}

\end{document}